\documentclass{amsart}
\usepackage{amssymb, amscd, amsmath, amsthm, epsf,graphicx , epsfig, latexsym}

\newtheorem{theorem}{Theorem}[section]

\newtheorem{definition}[theorem]{Definition}
 
\newtheorem{corollary}[theorem]{Corollary}

\newtheorem{conjecture}[theorem]{Conjecture}
\def\zz{{\bf Z}}
\def\qq{{\bf Q}}
\def\cc{{\bf C}}
\def\ff{{\bf F}}
\def\rr{{\bf R}}
\def\calc{\mathcal{C}}
\def\cala{\mathcal{A}}
 
\def\calg{\mathcal{G}} 
\def\calf{\mathcal{F}}

\def\co{\colon}

\newcommand{\fig}[2] { \includegraphics[scale=#1]{#2}  }

\begin{document}


\thispagestyle{empty}
\centerline{\large{{\bf Table of Contents for the Handbook of Knot
Theory}}}

\

\centerline{William W. Menasco and Morwen B. Thistlethwaite, Editors}

\
\vskip1in

\begin{enumerate}

\item Colin Adams, {\it Hyperbolic Knots}
\item Joan S. Birman and Tara Brendle, {\it Braids: A Survey}
\item John Etnyre  {\it Legendrian and Transversal Knots}
\item Greg Friedman, {\it Knot Spinning}
\item Jim Hoste, {\it The Enumeration and Classification of Knots and
Links}
\item Louis Kauffman,  {\it Knot Diagramitics}
\item Charles Livingston,  {\it A Survey of Classical Knot Concordance}
\item Lee Rudolph,  {\it Knot Theory of Complex Plane Curves}
\item Marty Scharlemann,  {\it Thin Position in the Theory of Classical Knots}
\item Jeff Weeks,  {\it Computation of Hyperbolic Structures in Knot
Theory}
\end{enumerate}
\vfill
\pagebreak


\title{A Survey of Classical Knot Concordance}

\author{Charles Livingston}

\address{Department of Mathematics, Indiana University, Bloomington, IN 47405}

\email{livingst@indiana.edu}
\keywords{Knot concordance}

\subjclass{Primary 57M25; Secondary 11E39}
 
\date{\today}


\maketitle

In 1926 Artin~\cite{a} described the construction of certain knotted 2--spheres in
$\rr^4$.  The intersection of each of these knots with the standard
$\rr^3 \subset
\rr^4$ is a nontrivial knot  in $\rr^3$. Thus a natural problem  is to identify  which
  knots can occur as such slices of knotted 2--spheres.   Initially it seemed possible that
every knot is such a {\em slice knot}  and it wasn't until the    early 1960s that
Murasugi~\cite{mu} and Fox and Milnor~\cite{fox1,fm}   succeeded at proving that some knots are
not slice.  

Slice knots can be used to define an equivalence relation   on the set of knots in
$S^3$: knots $K$ and $J$ are equivalent if $K \# -J$ is slice. With this equivalence the set of
knots becomes a group, the {\em concordance group} of knots. Much progress has been made in
studying slice knots and the concordance group, yet some of the most easily asked questions
remain untouched.  

There are two related theories of concordance, one in the smooth category and the other
topological.  Our focus will be on the smooth setting, though the distinctions and main results
in the topological setting will be included.  Related topics must be excluded, in particular
 the study of link concordance.   Our focus lies entirely in the classical setting; higher
dimensional concordance theory is only mentioned when needed to understand the classical setting.

 \section{Introduction}

Two smooth knots, $K_0$ and $K_1$, in $S^3$ are called {\it concordant} if there is a smooth
embedding of $ S^1 \times [0,1]$ into $S^3 \times [0,1]$ having boundary
  the knots $K_0$ and $-K_1$ in $S^3 \times \{0\}$ and  $S^3 \times
\{1\}$, respectively.  Concordance is an equivalence relation, and the set of equivalence
classes forms a countable abelian group,
$\calc$, under the operation induced by connected sum. A knot represents the trivial
element in this group if it is {\it slice}; that is, if it bounds an embedded disk in the
4--ball.

The concordance group was introduced in 1966 by Fox and Milnor in~\cite{fm}, though earlier work
on slice knots was already revealing aspects of its structure.   Fox~\cite{fox1}
described the use of the Alexander polynomial to prove that the figure eight knot  is of order
two in $\calc$ and Murasugi~\cite{mu} used the signature of a knot to obstruct the slicing of a
knot, thus showing that  the trefoil is of infinite order in $\calc$.  (These results, along with much of the introductory material, is presented in greater detail in the body of this article.) The application of abelian
knot invariants (those determined by the cohomology of abelian covers or, equivalently, by the
Seifert form) to concordance culminated in 1969 with Levine's classification of higher
dimensional knot concordance,~\cite{le1, le2}, which applied in the classical dimension to give
a surjective homomorphism,
$\phi \co  \calc
\to 
\zz^\infty
\oplus \zz_2^\infty \oplus \zz_4^\infty$.

 In 1975 Casson and Gordon~\cite{cg1, cg2} proved that Levine's homomorphism is not an
isomorphism, constructing nontrivial elements in the kernel, and Jiang~\cite{ji} expanded on
this to show that the kernel contains a subgroup isomorphic to
$\zz^\infty$.  Along these lines it was shown in~\cite{l5}  that the kernel also
contains a subgroup isomorphic to
$\zz_2^\infty$.   
The 1980s saw two significant developments in the study of concordance. The first was based on
Freedman's work~\cite{fr, fq} studying the structure of topological 4--manifolds.  One
consequence was that methods of Levine  and those of  Casson--Gordon   apply in the topological
locally flat category, rather than only in the smooth setting.  More significant, Freedman
proved that all knots with trivial Alexander polynomial  are in fact slice in the topological locally flat category.

The other important development concerns the application of differential geometric techniques to
the study of smooth 4--manifolds, beginning with the work of Donaldson~\cite{d, dk} and
including the introduction of Seiberg--Witten invariants and their application to symplectic
manifolds, the use of the  Thurston--Bennequin invariant~\cite{am, ru3}, and   recent work of
Ozsv\'ath and Szab\'o~\cite{os}.   This work quickly led to the construction of smooth knots of
Alexander polynomial one that are not smoothly slice, along with a much deeper understanding of
related issues, such as the 4--ball genus of knots.  Using these methods it has recently been shown that the results of~\cite{os} imply that the kernel of Levine's homomorphism contains a summand isomorphic to $\zz$ and thus contains elements that are not divisible~\cite{l8}.
References are too numerous to enumerate
here; a few will be included as applications are mentioned.

Recent work of Cochran, Orr and Teichner,~\cite{cot1, cot2}, has revealed a deeper structure to
the knot concordance group.  In that work a filtration of $\calc$ is defined: $$ \cdots
\calf_{2.0}
\subset\calf_{1.5} \subset \calf_{1} 
\subset \calf_{.5} \subset \calf_{0} \subset \calc.$$ It is shown that
$\calf_{0}$ corresponds to knots with trivial Arf invariant, $\calf_{.5}$ corresponds to knots
in the kernel of
$\phi$ and all knots in 
$\calf_{1.5}$ have vanishing Casson--Gordon invariants. Using von Neumann
$\eta$--invariants, it has been proved in~\cite{ct} that  each quotient is infinite. This work
places Levine's obstructions and those of Casson--Gordon in the context of an infinite sequence
of   obstructions, all of which   reveal a finer structure to
$\calc$.  

\vskip.1in
\noindent{\sl Outline} Section 2 is devoted to the basic definitions related to concordance and
algebraic concordance.  In Section 3 algebraic concordance invariants are presented, including
the description of  Levine's homomorphism. Sections 4 and 5 present Casson-Gordon invariants
and their application.  In Section 6 the consequences of Freedman's work on topological surgery
in dimension four are described.  Section 7 concerns the application of the results of Donaldson
and more recent differential geometric techniques to concordance.  In Section 8 the recent work
of Cochran, Orr, and Teichner, on the structure of the topological concordance group are
outlined.  Finally, Section 9 presents a few outstanding problems in the study of knot
concordance.

\vskip.1in

\noindent{\sl Acknowledgments}  This survey benefitted from the   suggestions
of many readers.  In particular,  conversations with Pat Gilmer were very helpful. Also, the
careful reading of Section 8 by Tim Cochran, Kent Orr and Peter Teichner  improved the
exposition there.  General references concerning concordance which were of great benefit to me as
I learned the subject include~\cite{go, rol}.


\section{Definitions}

We will work in the smooth setting.  In Section \ref{topcat} there will be a discussion of the
necessary modifications and main results that apply in the topological locally flat category. 
Knots are usually thought of as isotopy classes of embeddings of
$S^1$ into
$S^3$.  However, to simplify the discussion of orientation and symmetry issues,    it is
worthwhile to begin with the following precise definitions of  knots, slice knots,
 concordance and Seifert surfaces. 
\vskip.2in

\subsection{Knot Theory and Concordance.}

\begin{definition}
\hfill
\begin{enumerate} 
\item  A {\it knot}  is  an oriented diffeomorphism class of a pair of oriented manifolds,
$K = (\Sigma^3,
\Sigma^1)$, where
$\Sigma^n$ is diffeomorphic to the $n$--sphere. 

\item  A knot is called {\it slice} if there is a pair
$(B^4, D^2)$ with $\partial(B^4,D^2) = K$, where $B^4$ is the   4--ball  and
$D^2$ is a smoothly embedded 2--disk.  

\item Knots $K_1$ and $K_2$ are called {\it concordant} if $K_1\ \# -K_2$ is slice.  (Here $-K$
denotes the knot obtained by reversing the orientation of each element of the pair and connected
sum is defined in the standard way for oriented pairs.) The set of concordance classes is denoted
$\calc$.  

\item A Seifert surface  for a knot $K$ is an oriented surface $F$ embedded in
$S^3$ such that $K = (S^3, \partial F)$.
\end{enumerate}
\end{definition}
 
   The basic theorem in the subject is the following. 
 
\begin{theorem} The set of concordance classes of knots forms a countable abelian
group, also  denoted $\calc$, with its operation induced by connected sum and with the unknot
representing the identity.\end{theorem}

Related to the notion of slice knots there is the stronger condition of being a ribbon knot.

\begin{definition} A knot $K$ is called ribbon if it bounds an embedded disk $D$ in
$B^4$ for which the radial function on the ball restricts to be a smooth Morse function
with no local maxima in the interior of $D$.
\end{definition}

There is no corresponding group of ribbon concordance.  Casson observed that for every slice
knot $K$ there is a ribbon knot $J$ such that $K \# J$ is ribbon.  Hence, if any equivalence
relation identifies ribbon knots, it also identifies all slice knots. There is however a notion
of ribbon concordance, first studied in~\cite{go3}.

\subsection{Algebraic Concordance.}

An initial understanding of $\calc$ is obtained via the algebraic concordance group, defined by
Levine  in terms of Seifert pairings.

\begin{definition}  A Seifert pairing for a knot $K$ with Seifert surface $F$ is the bilinear
mapping $$V  \co H_1(F) \times H_1(F) \to
\zz$$ given $V(x,y) =
\mbox{lk}(x,i_*y)$, where $\mbox{lk}$ denotes the linking number and
$i_*$ is the map induced by the positive pushoff, $i \co F \to S^3 - F$.\end{definition}  

\noindent (Here
and throughout, homology groups will be taken with integer coefficients unless indicated
otherwise.)  A Seifert matrix is the matrix representation of the Seifert pairing  with respect
to some free generating set for
$H_1(F)$.

If the transpose pairing $V^\tau$ is defined by $V^\tau(x,y) = V(y,x)$ then $V- V^\tau$
represents the unimodular intersection form on $H_1(F)$.  Hence, in general we define an
abstract Seifert pairing on a finitely generated free $\zz$--module $M$ to be a bilinear form
$V \co M
\times M
\to \zz$ satisfying $V - V^\tau$ is unimodular.  (In order for this to make sense for the
trivial knot with Seifert surface $B^2$,   the Seifert form on the 0--dimensional
$\zz$--module is defined to be unimodular.)

\begin{definition} An abstract Seifert form $V$ on $M$ is metabolic if $M = M_1
\oplus M_2$ with $M_1 \cong M_2$ and $V(x,y) = 0$ for all $x$ and $y
\in M_1$.  Such an
$M_1$ is called a metabolizer for $V$.

\end{definition}

\begin{theorem}If $K$ is slice and $F$ is a Seifert surface for $K$, then the associated Seifert
form is metabolic.
\end{theorem}

\begin{proof}Let $D$ be a slice disk for $K$.  The union $F \cup D$ bounds a 3--manifold
$R$ embedded in $B^4$.  Such an $R$ can be constructed explicitly, or an obstruction theory
argument can be used to construct a smooth mapping $B^4 - D
\to S^1$ which has $F \cup D$ as the boundary of the pull-back of a regular value.  (Note
that this construction depends on the triviality of the normal bundle to $D$.)

A   duality argument implies that $\mbox{rank(ker}( H_1(F) \to H_1(R))) =
\frac{1}{2}\mbox{rank}(H_1(F))$.  For any $x$ and $y$ in that kernel,  $V(x,y) = 0$: since
$x$ bounds a 2--chain in $R$, $  i_*(x)$ bounds a 2--chain in $B^4 - R$   which is disjoint from
the chain bounded in $R$ by $y$. 

Since $V$ vanishes on this kernel, it vanishes on the summand $M$ generated by the kernel, and
hence $V$ is metabolic.

\end{proof}

\begin{corollary}\label{meta}If $K_1$ is concordant to $K_2$ and these knots have Seifert forms
$V_1$ and $V_2$ (with respect to arbitrary Seifert surfaces), then
$V_1 \oplus -V_2$ is metabolic.
\end{corollary}

In general, abstract Seifert forms $V_1$ and $V_2$ are called {\it algebraically concordant} if 
$V_1 \oplus -V_2$ is metabolic.  This is an equivalence relation.  (The proof is based on
cancellation: if
$V$ and $V\oplus W$ are metabolic, then so is $W$.  See~\cite{ke}.) 

\begin{theorem} The set of algebraic concordance classes forms a group, denoted
$\calg$, with its operation induced by direct sum.  The trivial
0--dimensional $\zz$--module serves as the identity.

\end{theorem} 

In the following theorem, defining Levine's homomorphism, we temporarily use the notation $[K]$
to denote concordance class of a knot and $[V_F]$ to represent the algebraic concordance class
of a Seifert form associated to an arbitrarily chosen Seifert surface $F$ for
$K$. 

\begin{theorem} The  map $\phi \co  \calc \to \calg$ defined by
$\phi([K]) = [V_F]$ is a surjective homomorphism.
\end{theorem}
 \begin{proof}That this map is well-defined follows from the previous discussion and, in
particular, Corollary~\ref{meta}.  Surjectivity follows from an explicit construction of a
surface with desired Seifert form~\cite{bz,se1}.
\end{proof}

\section{Algebraic Concordance Invariants}\label{aci}

In~\cite{le1} Levine defined a collection of homomorphisms from
$\calg$ to the groups
$\zz$, $\zz_2$ and $\zz_4$.  These can be properly combined to give an isomorphism
$\Phi$ from
$\calg$ to the infinite direct sum $\zz^\infty \oplus \zz_2^\infty
\oplus
\zz_4^\infty$.   The proof of this will be left to~\cite{le1}.

We should remark  that what Levine actually did was to classify the rational algebraic
concordance group, based on rational matrices.  He also showed that the integral group injects
into the rational group, with image sufficiently large to contain $\zz^\infty
\oplus \zz_2^\infty \oplus
\zz_4^\infty$.  Stolzfus, in~\cite{st}, completed the classification of the integral concordance
group.   

In this section we will describe a collection of invariants that are sufficient to show that
$\calg$ contains a summand isomorphic to 
$\zz^\infty \oplus
\zz_2^\infty \oplus
\zz_4^\infty$. The invariants  will be applied to a particular family of knots, which we now
describe.

\begin{figure}[h]
 \centerline{  \fig{.7}{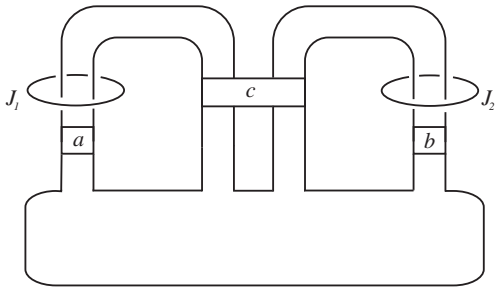} }  
 \caption{The knot $K(a,b,c)$.} \label{kabc}
 \end{figure}

\begin{figure}[h]
 \centerline{  \fig{.7}{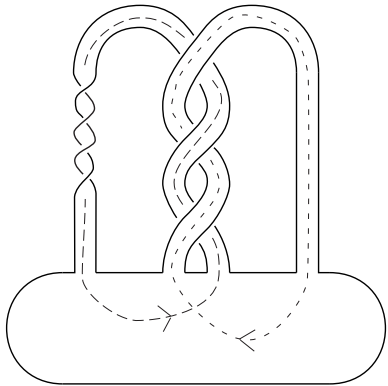} }  
 \caption{The knot $K(2,0,3)$.} \label{k203}
 \end{figure}

Figure~\ref{kabc} illustrates a basic knot that we denote
$K(a,b,c)$; the curves
$J_1$ and $J_2$ can be ignored for now.  The integers $a$ and $b$ indicate the number of full
twists in each band.  The integer $c$ is odd and represents the number of half twists between
the bands; those twists between the bands are so placed as to not add twisting to the individual
bands.   Figure~\ref{k203}  illustrates a particular example, $K(2,0,3)$, along with a basis for the first homology of the Seifert surface, indicated with oriented dashed curves on the surface.

 The knot $K(a,b,c)$  bounds a genus one Seifert surface  with Seifert form represented by the
following matrix with respect to the indicated basis of
$H_1(F)$.   
$$\left(  \begin{tabular}{c   c} 
     $a$ & $(c+1)/2$  \\ 
  $(c-1)/2$  & $b$ \\
\end{tabular}
\right) $$

\subsection{Integral Invariants, Signatures.} Let $V$ be a Seifert matrix and $V^\tau$ its
transpose. If
$\omega$ is a unit complex number that is not a root of the Alexander polynomial of $V$,
$\Delta_V(t) =
\det(V - tV^\tau)$, then the form $V_{\omega} =
\frac{1}{2}(1-\omega)V +\frac{1}{2} (1- \overline{\omega})V^\tau$ is nonsingular.  In this
case, if
$V$ is metabolic, the signature of
$V_\omega$ is 0.  To adjust for the possibility of the  $V_\omega$ being singular,  for general
$\omega$ on the unit circle  the signature
$\sigma_\omega(V)$ is defined to be the limiting average of the signatures of
$V_{\omega_+}$ and $V_{\omega_-}$, where $\omega_+$ and $\omega_-$ are unit complex numbers
approaching $\omega$ from different sides.  For all $\omega$,
$\sigma_\omega$ defines a homomorphism from $\calg$ to $\zz$.  It is onto $2\zz$ if
$\omega \ne 1$.  For the set of
$\omega$ given by roots of unity
$e^{2 \pi i /p}$ where $p$ is a prime, the functions $\sigma_\omega$ are  
  independent on $\calg$ (this can be seen using the $b$--twisted doubles of the unknot,
$K(1,b,1), b>0$), and hence together these give a map of
$\calg$ onto
$\zz^\infty$. In the case of $\omega = -1$, this signature, defined by Trotter~\cite{tr1},
was shown to be a concordance invariant by Murasugi~\cite{mu}.  The more general formulation is
credited to Levine and Tristram~\cite{tr} and is  referred to as the Levine-Tristram signature.

In~\cite{ka, kt1}  the identification of these signatures with signatures of the branched covers
of $B^4$  branched over a pushed in Seifert surface of a knot was made.  In~\cite{cl} it was
shown that the set of $\sigma_\omega$ over all $\omega$ with positive imaginary part are
independent.

\subsection{The Arf Invariant: $\zz_2$.}

Given a $(2g)\times (2g)$ Seifert matrix $V$ one defines a
$\zz_2$--valued quadratic form on $\zz_2^{2g}$ by
$q(x)= xVx^\tau$.  This is a nonsingular quadratic form in the sense that $q(x+y) - q(x)
-q(y) = x\cdot y$ where the nonsingular bilinear pairing $x \cdot y$ is given by the matrix
$V +V^\tau$.  (Recall that the determinant of $V+V^\tau$ is odd.)

The simplest definition  of the Arf invariant of a nonsingular quadratic form on a
$\zz_2$--vector space $W$ is that Arf$(q) = 0$ or Arf($q) = 1$ depending on whether $q$ takes
value 0 or 1, respectively, on a majority of elements in $W$.  See for instance~\cite{br}
or~\cite{hm}.  The Arf invariant defines a homomorphism on the Witt group of $\zz_2$ quadratic
forms and in particular vanishes on metabolic forms. Hence, the Arf invariant gives a well
defined
$\zz_2$--valued homomorphism from $\calg$ to $\zz_2$.

This invariant was first defined by Robertello in~\cite{ro}.  Murasugi~\cite{mu2} observed that
Arf($V$) = 0 if and  only if
$\Delta_V(-1) = \pm 1 \mod 8$.

\subsection{Polynomial Invariants: $\zz_2$.}
  The Alexander polynomial of a Seifert matrix is defined to be
$\Delta_V(t) =
\det(V - tV^\tau) \in \zz[t, t^{-1}]$. If different Seifert matrices associated to  the same
knot are used to compute an Alexander polynomial, the resulting polynomials will differ by
multiplication by a unit in
$\zz[t, t^{-1}]$, that is by $\pm t^n$ for some $n$.  Hence, two Alexander polynomials are
considered equivalent if they differ by multiplication   by $\pm t^n$ for some $n$.  

If $V$ is metabolic, then $\Delta_V(t) = \pm t^n f(t)f(t^{-1})$ for some integral polynomial
$f$.   For concordance considerations, if
$p(t)$ is an irreducible symmetric polynomial ($p(t^{-1}) = \pm t^n p(t)$) then the exponent of
$p(t)$ in the irreducible factorization of $\Delta_V(t)$ taken modulo 2 yields a $\zz_2$
invariant of
$\calg$.  Milnor and Fox~\cite{fm} used this to define a surjective homomorphism of $\calg$ to
$\zz_2^\infty$. The knots $K(a,-a,1)$ (see Figure~\ref{kabc})   are of order at most 2 in
$\calc$ since for each, $K(a,-a,1) = -K(a,-a,1)$.  On the other hand, these have distinct
irreducible Alexander polynomials if
$a>0$.  The existence of an infinite summand of
$\calg$ isomorphic to $\zz_2 ^ \infty$ follows.  Note that the knot $K(1,-1,1)$ is the figure eight knot. 

\subsection{W($\qq)$: $\zz_2 $ and $\zz_4$ Invariants.}  

  The matrix
$V + V^\tau$  defines an element in the Witt group of $\qq$, $W(\qq)$.  We will now summarize
the theory of this Witt group  and associated Witt groups of finite fields.  Details can be
found in~\cite{hm}.  Notice that the determinant of $V + V^\tau$ is odd; hence, in the
following discussion we restrict attention to odd primes
$p$.

Recall that the Witt group of an arbitrary field $\ff$ consists of finite dimensional
$\ff$--vector spaces with nonsingular symmetric forms and forms $W_1
$ and $W_2$ are equivalent if $W_1 \oplus - W_2$ is metabolic.  Addition is via direct sums.  

There is a surjective homomorphism
$\oplus
\partial_p \co   W(\qq)
\to
\oplus W(\ff_p )$.   Here $\ff_p$ is the field with
$p$ elements, and the direct sums are over the set of all primes. For $p$ odd, the group
W$(\ff_p)$ is isomorphic to either
$\zz_2$  or
$\zz_4$, depending on whether $p$ is 1 or 3 modulo 4. We next define $\partial_p$ and then
discuss the invariants of $W(\ff_p)$. (For completeness, we note here that the kernel of $ \oplus
\partial_p$ is $W(\zz)$ which is isomorphic to $\zz$ via the signature~\cite{hm}.)

\subsubsection{Reducing to Finite Fields.} 

There is a simple algorithm giving the map $\partial_p$.  A symmetric rational matrix
$A$ can be diagonalized using simultaneous row and column operations, and this form
decomposes as the direct sum of forms: $\oplus_{i=1}^n (a_i p^{\epsilon_i})$ where gcd($a_i,p)
=1$,
$\epsilon_i = 1$ for $i \le m$ and $\epsilon_i = 0$ for $m+1 \le i
\le n$.  The map
$\partial_p$ takes $A$ to the $\ff_p$ form represented by the direct sum
$\oplus_{i=1}^m (a_i)$.

\subsubsection{$W(\ff_p)${\rm :} $\zz_2$ and $\zz_4$ Invariants.} For $p$ odd, any form on a
finite dimensional
$\ff_p$--vector space can be diagonalized  with $\pm 1$ as the diagonal entries.  In the Witt
group  the form represented by the matrix $(1) \oplus(-1)$ is  trivial.  A little more work
shows that the   form $4(1)$ is Witt trivial: find elements $a$ and $b$ such that $a^2 + b^2 =
-1$ and consider the subspace spanned by
$(1,0,a,b)$ and
$(0,1,b,-a)$.  Hence,  $W(\ff_p)$ is generated by $(1)$, an  element of order 2 or 4.

In the case that $p \equiv 1 \mbox{\ modulo\ } 4$, $-1$ is a square.  It follows quickly that
$W(\ff_p) \cong \zz_2$.  On the other hand, in the case that $p
\cong 3
\mbox{\ modulo\ } 4$, $-1$ is not a square, and $W(F_p) \equiv
\zz_4$.

As a simple example, if one starts with the Seifert form for the knot $K(1,-5,1)$,
$$V = \left(  \begin{tabular}{c   c} 
     $1$ & $ 1$  \\ 
  $0$  & $-5$ \\
\end{tabular}
\right),\hskip.3in V+V^\tau = \left(  \begin{tabular}{c   c} 
     $2$ & $ 1$  \\ 
  $1$  & $-10$ \\
\end{tabular}
\right) .  $$  Diagonalizing over the rationals yields 
$$V = \left(  \begin{tabular}{c   c} 
     $2$ & $ 0$  \\ 
  $0$  & $-(2)(3)(7)$ \\
\end{tabular}
\right) .$$ With $p = 3$ this form  maps to the element $(-14)$ of
$W(\ff_3)$, which  is equivalent to the form $(1)$, a generator of order 4.  The same is true
working with $p = 7$. 

As a consequence of the next theorem we will see that this particular form $V$ is actually of
order  four in
$\calg$.

\subsection{Quadratic Polynomials.}
  A special case of a theorem of  Levine (Section 23 of~\cite{le2}) gives the following result,
which implies in particular that the form just described is of order 4 in $\calg$.  

\begin{theorem}\label{lequad} Suppose that $\Delta_V(t)$ is an irreducible quadratic.  Then $V$
is of finite order in the algebraic concordance group if and only if 
$\Delta_V(1)\Delta_V(-1) <0$. In this case  $V$ is of order 4 if 
$|\Delta_V(-1)| = p^a q$ for some prime $p$ congruent to 3 modulo 4,
$a$ odd, and
$p$ and
$q$ relatively prime; otherwise it is of order 2.
\end{theorem}

\subsection{Other Approaches to Algebraic Invariants.} There are alternative approaches to
algebraic obstructions to a knot being slice that do not depend on Seifert forms.  For instance,
Milnor~\cite{mi} described signature invariants based on his duality theorem for the infinite
cyclic cover of a knot complement. The equivalence of these signatures and those of Tristram and
Levine is proved in~\cite{ma}. There is also an interpretation of the algebraic concordance
group in term of  the Blanchfield pairing of the knot.

\section{Casson-Gordon Invariants}

In the case that $K$ is algebraically slice, Casson-Gordon invariants offer a further
obstruction to a knot being slice.  We follow the basic description of~\cite{cg1}.

\subsection{Definitions.}

We begin by reviewing  the {\sl linking form}   on torsion$(H_1(M))$ for an
oriented 3--manifold
$M$.  If $x$ and $y$ are curves representing torsion in the first homology, then lk($x,y)$ is
defined to be $(d \cap y)/n \in \qq / \zz$, where $d$ is a 2-chain with boundary  $nx$.
Intersections are defined via transverse intersections of chains, and of course one must check
that the value of the linking form is independent of the many choices in its definition.  For a
closed oriented 3-manifold the linking form is nonsingular in the sense that it induces an
isomorphism from
 torsion$(H_1(M))$ to hom(torsion$(H_1(M), \qq / \zz)$.

Such a symmetric pairing on a finite abelian group, $l \co H \times H \to \qq/\zz$, is
called {\sl metabolic} with {\sl metabolizer} $L$ if the linking form vanishes on $L \times L$
for some subgroup $L$ with
$|L|^2 = |H|$. 

Let $ M_q $ denote the $q$--fold branched cover of $S^3$ branched over a given knot
$K$, and let $\overline{M}_q$ denote $0$--surgery on $M_q$ along
$\tilde{K}$, where
$\tilde{K}$ is the lift of $K$ to $M_q$.  Here
$q$ will be a prime power.

   Let $x$ be an element of self-linking 0 in $H_1(M_q)$ and suppose that $x$ is of prime
power order, say $p$.  Linking with $x$ defines a homomorphism
$\chi_x \co H_1(M_q) \to
\zz_{p}$.  Furthermore, $\chi_x$ extends to give a $\zz_p$--valued character
 on
$H_1( \overline{M}_q )$ which vanishes on the meridian of
$\tilde{K}$.  In turn, this character extends to give
$\overline{\chi}_x \co  H_1( \overline{M}_q ) \to
\zz_p \oplus
\zz$.   Since $x$ has self-linking 0, bordism theory implies that the pair
$(\overline{M}_q,
\overline{\chi}_x)$ bounds a 4--manifold, character, pair, $(W,
\eta)$. 

More generally, for any character $\chi \co H_1(M_q) \to
\zz_{p}$, there is a corresponding character $\overline{\chi}  \co  H_1( \overline{M}_q ) \to
\zz_p \oplus
\zz$.  This character might not extend to a 4--manifold, but since the relevant bordism
groups are finite, for  some multiple $r\overline{M}_q $ the character given by
$\overline{\chi}$ on each component does extend to a 4--manifold, character pair, $(W,
\eta)$.

Let $Y$ denote the $\zz_p \times \zz$ cover of $ {W}$ corresponding to
$\eta$.  Using the action of  $\zz_p \times \zz$ on $H_2(Y,\cc)$  one can form the twisted
homology group $H_2^t( {W} , \cc) = H_2(Y,\cc)\otimes_{\cc[\zz_p \times
\zz]}\cc(t)$.  (The action of $\zz_p$ on
$\cc(t)$ is given by multiplication by $e^{2 \pi i /p}$.) There is a  nonsingular hermitian
form on $H_2^t(W,\cc)$ taking values in
$\cc(t)$. The Casson-Gordon invariant is defined to be the difference of this form and the
intersection form of
$H_2(W,\cc)$, both tensored with $\frac{1}{r}$, in $ W(\cc[t,t^{-1}]) \otimes \qq$.  (In
showing that this Witt class yields a  well-defined obstruction to slicing a knot, the  fact
that
$\Omega_4(\zz_p
\oplus
\zz) $ is nonzero appears, and as a consequence one must tensor with
$\qq$ to arrive at a well defined invariant, even in  the case of $\chi_x$ in which it is
possible to take $r =1$.)  

\begin{definition} The Casson-Gordon invariant $\tau (M_q,\chi)$ is the class
$(H_2^t( {W},\cc) -  H_2( {W}) , \cc))\otimes \frac{1}{r} \in W(\cc(t)) \otimes
\qq$.\end{definition}

\subsection{Main Theorem.}

The main theorem of~\cite{cg1} states:

\begin{theorem}If $K$ is slice, there is a metabolizer $L$ for the linking form on
$H_1(M_q)$ such that, for each prime power $p$ and each element $x
\in L$ of order $p$, $\tau(M_q, \chi_x)  = 0$.
\end{theorem}

The proof shows that if $K$ is slice with slice disk $D$, then covers of $B^4 - D$ can be used
as the manifold $W$, and for this
$W$ the invariant vanishes.

\vskip.1in 
\noindent{\bf Comment.}  There are a number of extensions of this theorem. With care the
definition of the Casson-Gordon invariant can be refined and
$\tau$ can be viewed as taking values in $W(\qq[\zeta_p](t)) \otimes
\zz[\frac{1}{p}]$. This
 yields finer invariants; see for instance~\cite{gl2}.  The observation that $L$ can be
assumed to be equivariant with respect to the deck transformation of $M_q$ can give stronger
constraints; see for example~\cite{kl}.    In~\cite{skim} it is demonstrated that a
factorization of the Alexander polynomial of a knot yields further constraints on the
metabolizer $L$.

\subsection{Invariants of $W(\cc(t)) \otimes \qq$.}

In the next section we will describe examples of algebraically slice knots which can be proved
to be nonslice using Casson-Gordon invariants.  We conclude this section with a description of
the types of algebraic invariants associated to the Witt group $W(\cc(t))
\otimes \qq$.
\subsubsection{Signatures.}  Let $\xi$ be a unit complex number.  Let $A \otimes
\frac{a}{b} 
\in   W(\cc(t)) \otimes \qq$. Then $A$ can be represented by a matrix of rational functions,
$A(t)$.  The signature $\sigma_\xi(A \otimes
\frac{a}{b})$ is defined, roughly, to be $\frac{a}{b}\sigma(A(\xi))$ where
$\sigma$ denotes the standard hermitian signature. There is the technical point arising that
$A(\xi)$ might be singular, so the precise definition of
$\sigma_\xi(A \otimes
\frac{a}{b})$ takes the two-sided average over unit complex numbers close to
$\xi$.  This limit is defined to be the {\em Casson-Gordon signature invariant},
$\sigma_{\xi}(K,\chi)$.  For $\xi =0$ this is abbreviated as $\sigma(K,\chi)$.

\subsubsection{Discriminants.}  If the matrix $A(t)$ represents $0
\in W(\cc(t))$, the discriminant, $\mbox{dis}(A(t)) =(-1)^{k}
\mbox{det}(A(t))$ (where $k$ is half the dimension of $A$) will be of the form
$f(t)\overline{f}(t)$ for some rational function
$f$.  Let $g(t) = t^2 + \lambda t +1, |\lambda| > 2$ be an irreducible real symmetric
polynomial.  It follows that for a matrix
$A(t)$, the exponent of $g(t)$ in the factorization of
$\mbox{dis}(A(t))$ gives a $\zz_2$--valued invariant of the Witt class of $A(t)$.  More
generally, in the case that $p$ is odd,  the exponent of $g$ in the determinant of $ A(t)^a$
gives a $\zz_2$ invariant of the class represented by $A(t)
\otimes \frac{a}{b}$ in $W(\cc(t) ) \otimes
\zz[\frac{1}{p}]$.  

These discriminants were first discussed in  unpublished work of Litherland~\cite{lit1}.  Later
developments and applications include~\cite{gl2}.  

In~\cite{kl3, kl} a three-dimensional approach to the definition of   Casson-Gordon discriminant
invariants is presented.  In short, the representation
$\overline{\chi} \co H_1(\overline{M}_q) \to \zz_p \times \zz$ determines the twisted homology
group: $H_1^t(\overline{M}_q,
\qq(\zeta_p)[t])$.  This is a
$\qq[\zeta_p][t]$ module, and the discriminant of the Casson-Gordon invariant is
 given by the order of this module. Although this three-dimensional approach
  does not give the signature invariant, it has the advantage of being completely algorithmic in
computation, via a  procedure first developed in~\cite{lin, wa} and applied in~\cite{kl3, kl, 
klv}.  A computer implementation of that algorithm facilitated the  classification of the order
of low-crossing number knots in concordance~\cite{tam} and the proof that most low-crossing
number knots which are not reversible are not concordant to their reverses, in~\cite{tam2}.

In a different direction, we note that some effort has been made in removing the restriction on
prime power covers and characters.  In the case of ribbon knots, it was known that stronger
results could be attained.  Recent work of Taehee Kim~\cite{kim} has developed examples of
nonslice algebraically slice knots for which all prime power branched covers are homology
spheres.  Other work in this realm includes that of Letsche~\cite{let} and recent   work of
Friedl~\cite{fri1, fri2}.

\section{Companionship and Casson-Gordon Invariants}

In Casson and Gordon's original work the computation of Casson-Gordon invariants was quite
difficult, largely limited to restricted classes of knots.   
Litherland~\cite{lit2}  studied the behavior of these invariants under companionship and,
independently, Gilmer~\cite{gi} found interpretations of particular Casson-Gordon invariants in
terms of signatures of simple closed curves on a Seifert surface for a knot.   Further work
addressing companionship and Casson-Gordon invariants includes~\cite{ab}.  In this section we
describe the general theory and its application to genus one knots.

\subsection{Construction of Companions.} Let $U$ be an unknotted    circle in the complement of
a knot $K$.  If $S^3$ is modified by removing a neighborhood of $U$ and replacing it with the
complement of a knot $J$ in $S^3$ (via a homeomorphism of boundaries that identifies the
meridian of $J$ with the longitude of $U$ and vice versa) then the resulting manifold is again
diffeomorphic to $S^3$.  The image of $K$ in this manifold will be denoted $K(J)$ (the choice
of  $U$ will be suppressed in the notation).  In the language of classical knot theory,
$K(J)$ is a satellite knot with companion $J$ and satellite $K$.

If $M_q$ is the $q$--fold branched cover of $S^3$ branched over $K$, then $U$ has
$q' $ lifts, denoted $U_i$, $i = 1, \ldots , q'$, where $q' = \gcd (q, \mbox{lk}(U, K))$.  It
follows that
$M'_q$, the
$q$--fold branched cover of $S^3$ branched over $K(J)$, is formed from $M_q$ by removing
neighborhoods of the $U_i$ and replacing each with the $q/q'$--cyclic cover of the complement of
$J$.
 If
$\chi
$ is a
$\zz_p$--valued homomorphism on $H_1(M_q)$, there is a naturally associated homomorphism
$\chi'$ on $H_1(M'_q)$. 

\subsection{Casson-Gordon Invariants and Companions.}

 In the case that lk$(U, K) = 0$, we have the following theorem of Litherland~\cite{lit2}.  (See
also~\cite{gl1}.)

 \begin{theorem}\label{companion}  In the situation just described, with lk$(U,K) = 0$,
$$\sigma(K(J),
\chi') = 
 \sigma(K,\chi) + \sum_{i=1}^q \sigma_{\chi(U_i)/p} (J).$$
 \end{theorem}

The main idea of the proof is fairly simple.  If $({W}, \eta)$ is the chosen pair bounding
$(\overline{M}_q,
\overline{\chi})$ in the definition of the Casson--Gordon invariant, then for the new knot $K'$
a 4--manifold $W'$ can be built from $W$ by attaching copies of a 4--manifold with character
$(Y,\eta)$ bounding $0$--surgery on $J$ with its canonical representation to
$\zz$. Signatures of cyclic covers of $Y$ are related to the signatures of $J$.
A similar analysis can be done for   the discriminant of the Casson--Gordon invariant.  This was
detailed in~\cite{gl2}, and further explored in~\cite{klv} where it was no longer assumed that
$J$ was null-homologous.

\vskip.1in
\noindent{\bf Example.} Consider the knot illustrated in Figure~\ref{kabc} with
$a=0$,
$b=0$ and $c=3$.  The bands have knots $J_1$ and $J_2$ tied in them.  (We will also refer to the
pair of unknotted circles as $J_1$ and $J_2$ in this  situation, as the meaning is
unambiguous.)  Call the resulting knot $K(J_1, J_2)$.  The homology of the 2--fold cover is
isomorphic to
$\zz_3
\oplus
\zz_3$ with the linking form vanishing on the two summands.  Call generators of the summands
$x_1$ and $x_2$.  An analysis of the cover shows that $\chi_{x_1}$ is a
$\zz_3$--valued character that vanishes on the lifts of $J_1$ and takes value $\pm 1$ on the two
lifts of $J_2$.  Similarly for
$\chi_{x_2}$.

Since $K$ is slice, by the Casson--Gordon theorem, either
$\sigma(K,\chi_{x_1})$ or 
$\sigma(K,\chi_{x_1})$ must vanish.  Hence, using Theorem~\ref{companion}, if
$K(J_1,J_2)$ is slice, either $2\sigma_{1/3}(J_1)$ or
$2\sigma_{1/3}(J_2)$ must vanish.  By choosing $J_1$ and $J_2$ so that this is not the case, one
constructs basic examples of algebraically slice knots which are not slice.

\subsection{Genus One Knots and the Seifert Form.}

  Gilmer observed in~\cite{gi,gi2} that for genus one knots the computation of Casson-Gordon
invariants is greatly simplified.  Roughly, he interpreted the Casson-Gordon signature
invariants of an algebraically slice genus one knot in terms of the signatures of knots tied in
the bands of the Seifert surface.  The previous example offers an illustration of the appearance
of these signatures.  This work is now most easily understood via the use of companionship just
described.

In short, if an algebraically slice knot $K$ bounds a genus one Seifert surface
$F$, then some nontrivial primitive class in $H_1(F)$ has trivial self-linking with respect to
the Seifert form.  If that class is represented by a curve $\alpha$, the surface can be deformed
to be a disk with two bands attached, one of which is tied into the knot
$\alpha$.  If a new knot is formed by adding the knot $-\alpha$ to the band, the knot becomes
slice  and certain of its Casson-Gordon invariants will vanish.  However, the previous results on
companionship determine  how the modification of the knot changes the Casson-Gordon invariant. 
The situation is made somewhat more delicate in that
$\alpha$ is not unique: for genus one algebraically slice knots there are two metabolizers.  The
following represents the sort of result that can be proved.  

\begin{theorem} Let $K$ be a genus one   slice knot.  The Alexander polynomial of
$K$ is given  
$(a t - (a+1))((a+1)t - a)$ for some $a$.   For some simple closed curve
$\alpha$  representing a generator of  a metabolizer  of the Seifert form and for some infinite
set of primes powers $q$, one has
$$\sum_{i = 1}^q \sigma_{bm^i/p}(\alpha) = 0$$  for all prime power divisors $p$ of
$ (a-1)^q - a^q$, and for all integers $b$. 

\end{theorem}

(The appearance of the term $(a-1)^q - a^q$ represents the square root of the order of the
homology of the $q$--fold branched cover.) Since the sum is taken over a coset of the
multiplicative subgroup of
$\zz_p$, by combining these cosets one has the following.

\begin{corollary}If $K$ is a genus one slice knot with nontrivial Alexander polynomial, then for
some simple closed curve
$\alpha$ representing a generator of a metabolizer of the Seifert form, there is an infinite set
of prime powers
$p$ for which 
$$\sum_{i = 1}^{p-1} \sigma_{i/p}(\alpha) = 0.$$ 
\end{corollary}

A  theorem of Cooper~\cite{co}  follows quickly:

\begin{corollary}If $K$ is a genus one slice knot with nontrivial Alexander polynomial, then for
some simple closed curve
$x$ representing a generator of a metabolizer of the Seifert form,  
$$\int_0^{1/2}   \sigma_{t}(x) dt = 0.$$ 
\end{corollary}

\noindent (This theorem reappears in~\cite{cot2} where the integral is reinterpreted as a
metabelian von Neumann signature of the original knot $K$, giving a direct reason why it is a
concordance invariant.  For more on this, see Section~\ref{filtrat}.)

\vskip.1in 

\noindent{\bf Example.}  Consider the knot $K(0,0,3)$, as in Figure~\ref{kabc}.   Replacing the
curves labeled $J_1$ and $J_2$ with the complements of knots $J_1$ and
$J_2$ yields a knot for which the metabolizers of the Seifert form are represented by the knots
$J_1$ and $J_2$.  The knot is algebraically slice, but by the previous corollary, if both of the
knots have  signature functions with nontrivial integral, the knot is not slice.


\section{The Topological Category}\label{topcat}

In~\cite{fr} Freedman developed surgery theory in the category of topological 4--manifolds,
proving roughly that  for manifolds with fundamental groups that are not too complicated (in
particular, finitely generated abelian groups) the general theory of higher-dimensional surgery
descends to dimension 4.  The most notable consequence of this work was the proof the  
4--dimensional Poincar\'e Conjecture: a closed topological 4--manifold that is homotopy
equivalent to the 4--sphere is homeomorphic to the 4--sphere. 

Two significant contributions to the study of concordance quickly followed from Freedman's
original paper.  The first of these, proved in~\cite{fq}, is that a locally flat surface in a
topological 4--manifold has an embedded normal bundle.  The use of such a normal bundle was
implicit in the proof that slice knots are algebraically slice.  It is also used in a key step
in the proof of the Casson--Gordon theorem, as follows.  Casson--Gordon invariants of slice
knots are shown to vanish via the observation that for a slice knot $K$, if 0--surgery is
performed on $K$, the resulting 3--manifold $M(K,0)$ bounds a homology $S^1 \times B^3$, $W$. 
This
$W$ is constructed  by removing a tubular neighborhood of a slice disk for $K$ in the
$4$--ball.  The existence of the tubular neighborhood is equivalent to the existence of the
normal bundle.

In a different direction, Freedman's theorem implied that in the topological locally flat 
category all knots of Alexander polynomial one are slice.   To  understand why this is
a consequence, note first the following.

\begin{theorem}For a knot $K$, if $M(K,0)$ bounds a homology $S^1
\times B^3$,
$W$, with $\pi_1(W) = \zz$ then $K$ is slice.
\end{theorem} 

\begin{proof} We have that  $M(K,0)$ is formed from $S^3$ by removing a solid torus and
replacing it with another solid torus.  Performing 0--surgery on the core, $C$, of that solid
torus returns
$S^3$.  Attach a 2--handle to
$W$ with framing 0 to $C$.  The resulting manifold is a homotopy ball with boundary
$S^3$, and hence, by the Poincar\'e conjecture, is homeomorphic to
$B^4$.  The cocore of that added 2--handle is a slice disk for the boundary of the cocore, which
can seen to be the original $K$.
\end{proof}

Freedman observed that a surgery obstruction to finding such a manifold $W$ is determined by the
Seifert form, and for a knot of Alexander polynomial one that is the only obstruction, and it
vanishes.  

\subsection{Extensions.}

Is it possible that more delicate arguments using 4--dimensional surgery might yield stronger
results, showing that other easily identified classes of algebraically slice knots are  
slice, based only on the Seifert form of the knot?  The following result indicates that the
answer is no.

\begin{theorem}If $\Delta_K(t)$ is nontrivial then there are two nonconcordant knots having that
Alexander polynomial.
\end{theorem}

This result was first proved in~\cite{l2} where there was the added constraint that the
Alexander polynomial  is not the product of cyclotomic polynomials $\phi_n(t)$ with $n$
divisible by three distinct primes.  The condition on Alexander polynomials is technical,
assuring that some prime power branched cover is not a homology sphere.  Taehee Kim~\cite{kim}
has shown this  condition is not  essential in particular cases,  and in unpublished work he has
shown that the result applies for all nontrivial Alexander polynomials.


\section{Smooth Knot Concordance}

In 1983 Donaldson~\cite{d}  discovered  new constraints on the intersection forms of smooth
4--manifolds. This and subsequent work soon yielded the following theorem.

\begin{theorem}\label{don} Suppose that $X$ is a smooth closed 4--manifold and
$H_1(X,
\zz_2) = 0$.  If the intersection form on $H_2(X)$ is positive definite then the form is
diagonalizable.  If the intersection form is even and definite, and hence of the type $nE_8
\oplus mH$, where
$H$ is the standard 2--dimensional hyperbolic form, then if
$n>0$, it follows that $m>2$.
\end{theorem}

This result is sufficient to prove that many knots of Alexander polynomial one are not slice. 
The details of any particular example cannot be presented here, but the connections with
Theorem~\ref{don} are easily explained.

Let $M(K,1)$ denote the 3--manifold constructed as 1--surgery on
$K$.  Then
$M(K,1)$ bounds the 4--manifold $W$ constructed by adding a 2--handle to the 4--ball along $K$
with framing 1.  If $K$ is slice, the generator of $H_2(W)$ is represented by a 2--sphere with
self-intersection number 1.  A tubular neighborhood of that sphere can be removed and replaced
with a 4--ball, showing that
$M(K,1)$ bounds a homology ball, $X$.  If  
$M(K,1)$ also bounds a 4--manifold $Y$ (say simply connected) with intersection form of the type
obstructed by Theorem~\ref{don}, then a contradiction is achieved using the union of $X$ and
$Y$.  

As an alternative approach, notice that if $K$ is slice, the 2--fold branched cover of
$S^3$ branched over $K$, $M_2$, bounds the
$\zz_2$--homology ball formed as the 2--fold branched cover of
$B^4$ branched over the slice disk.  Hence, if $M_2$ is known to bound a  simply connected 
4--manifold with one of the forbidden forms of Theorem~\ref{don}, then again a contradiction is
achieved.  

It seems that prior to Donaldson's work it was known that either of these approaches would be
applicable to proving that particular polynomial one knots are not slice, but these arguments
were not published.  In particular, following the announcement of Donaldson's theorem it
immediately was known that  the pretzel knot $K(-3,5,7)$ and the untwisted double of the trefoil
(Akbulut) are not slice.   Early papers  presenting details of such arguments include~\cite{gom}
where it was shown that there are topologically slice knots of infinite order in smooth
concordance.  See~\cite{cg} for further examples.

\vskip.1in

\subsection{Further Advances.}

Continued advances in smooth 4--manifold theory have led to further understanding of the knot
slicing problem.  In particular, proving that large classes of Alexander polynomial one knots
are not slice has fallen to algorithmic procedures.  Notable among this work is that of
Rudolph~\cite{ru1, ru2, ru3}.  Here we outline briefly the approach using Thurston-Bennequin
numbers, as described by Akbulut and Matveyev in the paper~\cite{am}.

The 4--ball has a natural complex structure.  If a 2--handle is added to the 4--ball along a
knot $K$ with appropriate framing, which we call $f$ for now, the resulting manifold $W$ will
itself be complex.  According to~\cite{lm},
$W$ will  then embed in a closed Kahler manifold $X$.  Further restrictions on the structure of
$X$ are known to hold, and with these constraints the adjunction formula of Kronheimer and
Mrowka~\cite{km1,km2}  applies to show that no essential 2--sphere in $X$ can have
self-intersection greater than or equal to
$-1$.

On the other hand, if $K$ were slice and the framing $f$ of $K$ were greater than
$-2$, such a sphere would exist.  The appropriate framing  $f$ mentioned above depends on the
choice of representative of $K$, not just its isotopy class. If the representative is ${\bf
K}$,  then $f =
 tb({\bf K}) -1$, where $tb({\bf K}) $ is the Thurston-Bennequin number, easily computed from a
diagram for ${\bf K}$.  

Applying this, both Akbulut-Matveyev~\cite{am} and Rudolph~\cite{ru3} have given simple proofs
that, for instance, all iterated positive twisted doubles of the right handed trefoil are not
slice.

Although these powerful techniques have revealed a far greater complexity to the concordance
group than had been expected, as of yet they seem incapable of addressing some of the basic
questions: for instance the slice implies ribbon conjecture and problems related to torsion in
the concordance group.

\section{Higher Order Obstructions and the Filtration of $\calc$}\label{filtrat}

Recent work of Cochran, Orr, and Teichner  has demonstrated a deep structure to the topological
concordance group.  This  is revealed in a filtration of the concordance group by an infinite
sequence of subgroups:

$$ \cdots \calf_{2.0} \subset\calf_{1.5} \subset \calf_{1} 
\subset \calf_{.5}   \subset \calf_0 \subset \calc.$$
 This approach has successfully placed known obstructions to the slicing problem---the Arf
invariant, algebraic sliceness, and Casson-Gordon invariants---as the first in an infinite
sequence of invariants.    Of special significance is that each level of the induced filtration
of the concordance group has both an algebraic interpretation and a geometric one.      Here we
can offer   a simplified view of the motivations and consequences of their work, and in that
interest will focus   on the $\calf_n$ with $n$ a nonnegative integer.

To begin, suppose that $M(K,0)$,  0--surgery on a knot $K$, bounds a 4--manifold
$W$ with the homology type and intersection form of $S^1
\times B^3
\#_n S^2
\times S^2$. Such a $W$ will exist if and only if the Arf invariant of
$K$ is trivial.  Constructing one such  
$W$ is fairly simple in this case. Push a Seifert surface $F$ for
$K$ into
$B^4$ and perform surgery on
$B^4$ along a set of curves on $F$ representing a basis of a metabolizer for its intersection 
form, with the additional condition that it represents a metabolizer for the $\zz_2$--Seifert
form. (Finding such a basis is where the Arf invariant condition appears.)  When performing the
surgery, the surface
$F$ can be  ambiently surgered to become a disk, and the complement of that disk is the desired
$W$.   

 If a generating set of a metabolizer for the intersection form on
$H_2(W)$ could be represented by disjoint embedded 2--spheres, then surgery could be performed
on $W$ to convert it into a homology $S^1
\times B^3$.  It would quickly follow that $K$ would   be slice in a homology 4--ball bounded by
$S^3$.  

In the higher dimensional analog (of the concordance group of knotted
$(2k-1)$--spheres in $S^{2k+1}$, $k >1$),   there  is   
  an obstruction (to finding this family of spheres)  related to the twisted intersection form on
$H_{k+1}(W,\zz[\pi_1(W)])$, or, equivalently, related to the intersection form on the universal
cover of $W$.  In short, the intersection form of
$W$ should have a metabolizer that lifts to a metabolizer in the universal cover of $W$. In this
higher dimensional setting, if the obstruction vanishes then, via the Whitney trick, the
metabolizer for
$W$ can be realized by embedded spheres and $W$ can be surgered as desired.  This viewpoint on
knot concordance has its roots in the work of Cappell and Shaneson~\cite{cs}.

Whether in   high dimensions or in the classical   setting, the explicit construction of
a
$W$ described earlier in this section  yields a
$W$ with cyclic fundamental group.  This obstruction is thus determined solely by the infinite
cyclic cover and vanishes for algebraically slice knots.  Of course, in higher dimensions
algebraically slice knots are slice.  Clearly something more is needed in the classical case.

In light of the Casson-Freedman  approach to 4--dimensional surgery theory, in addition to
finding   immersed spheres representing a metabolizer for $W$, one needs to find appropriate
dual spheres in order to convert the  immersed spheres into embeddings. The Cochran-Orr-Teichner
filtration can be interpreted as a sequence of obstructions to finding a family of spheres and
dual spheres. To describe
 the filtration, we denote 
$\pi^{(0)} =
\pi =
\pi_1(W)$ and let $\pi^{(n)} $ be the derived subgroup:
$\pi^{(n+1)} = [\pi^{(n)}, \pi^{(n)}]$.

\begin{definition} A knot $K$ is called $n$--solvable if there exists a (spin) 4--manifold
$W$ with boundary $M(K,0)$ such that: (a) the inclusion map
$H_1(M(K,0)) \to H_1(W)$ is an isomorphism; (b) the intersection form on $H_2(W,\zz[\pi /
\pi^{(n)}])$ has a dual pair of self-annihilating submodules (with respect to intersections and
self-intersections), $L_1$ and $L_2$; and (c) the images of $L_1$ and
$L_2$ in $H_2(W)$ generate $H_2(W)$. 

\end{definition}

\noindent (Here and in what follows we leave the description of
$n.5$--solvability to~\cite{cot1}.)

There are the  following basic corollaries of the work in~\cite{cot1}.

\begin{theorem} If the Arf invariant of a knot $K$ is 0, then $K$ is
$0$--solvable.  If $K$ is $1$--solvable, $K$ is algebraically slice. If
$K$ is $2$--solvable, Casson--Gordon type obstructions to $K$ being slice vanish. If $K$ is
slice, $K$ is $n$--solvable for all $n$.
\end{theorem}

One of the beautiful aspects of~\cite{cot1} is that this very algebraic formulation is closely
related to the underlying topology. For those familiar with the language of Whitney towers and
gropes, we have the following theorem from~\cite{cot1}.
\begin{theorem}If $K$ bounds either a Whitney tower or a grope of height
$n+2$ in $B^4$, then $K$ is $n$--solvable.
\end{theorem}

Define $\calf_n$ to be the subgroup of the concordance group consisting of
$n$--solvable knots.  One has the filtration (where  we have dropped the
$n.5$--subgroups).

 $$ \cdots  \calf_{3}  \subset\calf_{2}  \subset \calf_{1}  \subset \calf_{0}
  \subset \calc.$$ Beginning with~\cite{cot1} and culminating in~\cite{ct}, there is the
following result.

\begin{theorem} For all $n$, the quotient group $\calf_n /
\calf_{n+1}$ is infinite and
$\calf_2 / \calf_3$ is infinitely generated.
\end{theorem}  

Describing the invariants that provide obstructions to a knot being in
$\calf_{n}$ is beyond  the scope of this survey.  However, two important aspects should be
mentioned.  First,~\cite{cot1} identifies a connection between
$n$-solvability and the structure and existence of metabolizers for   linking forms on
$$H_1(M(K,0),
  \zz[\pi_1(M(K,0)) /  \pi_1(M(K,0))^{ (k)}]) ,\hskip.1in  k\le n,$$ generalizing the fact that
for algebraically slice knots the Blanchfield pairing of the knot vanishes. 

The second aspect of proving the nontriviality of $\calf_n /
\calf_{n+1}$ is the appearance of von Neumann signatures for solvable quotients of the knot
group. Though difficult to compute in general,~\cite{cot1} demonstrates that if $K$ is built as
a satellite knot, then in special cases, as with the Casson-Gordon invariant,   
 the value of this complicated invariant is related to the Tristram-Levine signature function
of the companion knot.  More precisely,  
if a knot
$K$ is built from another knot by removing an unknot
$U$ that lies in
$\pi^{(n)}$ of the complement and replacing it with the complement of a knot $J$, then the change
in a particular    von Neumann
$\eta$--invariant of the
$\pi^{(n)}$--cover is related to the integral of the Tristram-Levine signature function of
$J$, taken over the entire circle. The Cheeger-Gromov estimate for these
$\eta$--invariants can then be applied to show   the nonvanishing of the invariant by
choosing
$J$ in a way that the latter integral exceeds the estimate.  This construction generalizes in a
number of ways the one used in applications of the Casson-Gordon invariant described earlier,
which applied only in the case that $U \in \pi^{(1)}$ and $U \notin \pi^{(2)}$.  Furthermore, the
Casson-Gordon invariant is based on a finite dimensional representation where here the
representation becomes infinite dimensional.  In the construction of \cite{ct} it is also
required that one work with a family of unknots; a single curve $U$ will not suffice.

\section{Three-dimensional Knot Properties and Concordance.}\label{3dprops}
\subsection{Primeness.} The first result  of the sort to be discussed here is the theorem of Kirby
and Lickorish~\cite{klk}:

\begin{theorem} Every knot is concordant to a prime knot.\end{theorem}

Shorter proofs of this were given in~\cite{l, na2}.  In these constructions it was shown that
the concordance can be chosen so that the Seifert form, and hence the algebraic invariants, of
the knot are unchanged.  Myers~\cite{my} proved that every knot is concordant to a  knot with
hyperbolic complement, and hence to one with no incompressible tori in its complement.  Later,
Soma~\cite{so} extended Myers's result  by showing that fibered knots are (fibered) concordant to
fibered hyperbolic knots.

In the reverse direction, one might ask if every knot is concordant to a composite knot, but the
answer here is obviously yes:   $K$ is concordant to $K \# J$, for any  slice knot $J$. 
However, when the Seifert form is taken into consideration the question becomes more
interesting.  Here we have the following example, the proof of which is contained
in~\cite[version 1]{l6}.

\begin{theorem}There exists a knot $K$ with Seifert form $V_K = V_{J_1}
\oplus V_{J_2}$, but $K$ is not concordant to a connected sum of knots with Seifert forms
$V_{J_1}$ and $V_{J_2}$.

\end{theorem}

Notice that by Levine's classification of higher dimensional concordance, such examples cannot
exist in dimensions greater than 3.

\subsection{Knot Symmetry: Amphicheirality.}

For the moment,  view a knot $K$ formally  as a smooth oriented pair $(S, K)$ where $S$ is
diffeomorphic to $S^3$ and $K$ is diffeomorphic to
$S^1$.  Equivalence is up to orientation preserving diffeomorphism.  (In dimension three it does
not matter whether the smooth or locally flat topological category is used.)

\begin{definition} A knot $(S,K)$ is called reversible (or invertible), negative amphicheiral,
or positive amphicheiral, if it is equivalent to
$K^r = (S, {-}K)$, $-K = (-S, {-}K)$, or $-K^r = (-S,K)$, respectively.  It is called strongly
reversible, strongly positive amphicheiral, or strongly negative amphicheiral if there is an
equivalence that is an involution. 
\end{definition}   

Each of these properties constrains the algebraic invariants of a knot, and hence can constrain
the concordance class of a knot.  For instance, according to Hartley~\cite{h}, if a knot $K$ is
negative amphicheiral, then its Alexander polynomial satisfies
$\Delta_K(t^2) = F(t)F(t^{-1})$ for some symmetric polynomial  $F$.  It follows quickly from the
condition that slice knots have polynomials that factor as $g(t)g(t^{-1})$ that if a knot $K$ is
concordant to a negative amphicheiral knot, $\Delta_K(t^2)$ must factor as
$F(t)F(t^{-1})$.  Further discussion of amphicheirality and knot concordance is included
in~\cite{cm}, where the focus is on higher dimensions, but some results apply in dimension three.

  \vskip.1in

\noindent{\bf Example.} Let $K$ be a knot with Seifert form $$V_a =
 \left(  \begin{tabular}{c   c} 
      1  &   1   \\ 
   0   &  $-a$  \\
\end{tabular}
\right).$$  If $a$ is positive, it follows from Levine's characterization of knots with
quadratic Alexander polynomial (Theorem~\ref{lequad}) that $K$ is of order two in the algebraic
concordance group if every prime of odd exponent in $4a +1$ is congruent to 1 modulo 4.  It
follows as one example  that any knot with  Seifert form $V_3$, for instance the 3--twisted
double of the unknot,  is of order 2 in algebraic concordance but is not concordant to a
negative amphicheiral knot. 

This gives insight into the following conjecture, based on a long standing question of
Gordon~\cite{go2}: 

\begin{conjecture}\label{gconject} If $K$ is of order two in
$\calc$, then
$K$ is concordant to a negative amphicheiral knot.
\end{conjecture}
 \noindent (Gordon's original question did not have the ``negative'' constraint in its
statement.) 

In a different direction, it was noted by Long~\cite{lo} that the example of a knot
$K$ for which $K \# -K^r$ is not slice (described in the next subsection) yields an example of a
nonslice strongly positive amphicheiral knot.  Flapan~\cite{fl} subsequently found a prime
example of this type.  It has since been shown that the concordance group contains   infinitely
many linearly independent such knots~\cite{l4}. 

\subsection{Reversibility and Mutation.}

Every knot is algebraically concordant to its reverse.  A stronger result, but the only proof in
print, follows from Long~\cite{lo}: if
$K$ is strongly positive amphicheiral then it is algebraically slice.  For any knot, $K
\# -K^r$ is strongly positive amphicheiral, so $K$ and $K^r$ are algebraically concordant. It is
proved in~\cite{li2} that  there are knots that are not concordant to their reverses.  Further
examples have been developed in~\cite{kl3, n3, tam2}.

Kearton~\cite{ke4} observed that since $K \# -K^r$ is a (negative) mutant of the slice knot $K \#
-K$, an example of a knot which is not concordant to its  reverse  yield an example  of mutation
changing the concordance class of a knot.  Similar examples for positive mutants proved harder
to find and were developed in~\cite{kl3, klv}.

\subsection{Periodicity.} A knot $K$ is called periodic if it is invariant under a periodic
transformation $T$ of $S^3$ with the fixed point set of $T$ a circle disjoint from $K$.  Some of
the strongest results concerning periodicity are those of Murasugi~\cite{mu3} constraining the
Alexander polynomials of such knots.  Naik~\cite{n1} used Casson--Gordon invariants to obstruct
periodicity for knots for which all algebraic invariants coincided with those of a periodic knot.

A theory of periodic concordance has been developed.  Basic results in the subject include those
of~\cite{ck, n2} obstructing knots from being periodically slice and those of~\cite{dn}
giving a characterization of the Alexander polynomials of periodically ribbon knots.  

\subsection{Genus.} The 4--ball genus of a knot $K$, $g_4(K)$, is the minimal genus of an
embedded surface bounded by $K$ in the 4--ball.  It is a concordance invariant of a knot
which is clearly bounded by its 3--sphere genus.  

This invariant has been studied extensively. It is known
to be bounded below by half the classical signature and the
Tristram-Levine signature~\cite{  kt1, ka, mu, tr}.  In the case that a knot is algebraically
slice, Gilmer developed bounds on the 4--ball genus using Casson-Gordon invariants~\cite{gi3}. 
In~\cite{kimliv} it is shown that for any pair of nonnegative integers $m$ and $n$ there is a
knot $K$ with a mutant $K^*$ such that $g_4(K) = m$ and $g_4(K^*) = n$; a knot and its
mutant are algebraically concordant. 
Beyond that, there are many results giving bounds on the 4-ball genus in the smooth setting
based on differential geometric results.  See for instance,~\cite{ru3, ta}. 

 Nakanishi~\cite{na} and Casson observed that there  are knots that bound surfaces of
genus one in the 4--ball but which are not concordant to knots of 3--sphere genus 1. 
In~\cite{l7} this observation was the starting point of the definition of the concordance genus
of a knot
$K$: the minimum genus among all knots concordant to $K$.  It is shown that this invariant can
be arbitrarily large, even for knots of 4--ball genus 1, and even among algebraically slice
knots.

\subsection{Fibering.}  A knot is called fibered if its complement is a surface bundle over
$S^1$.   It is relatively easy to see that not all knots are concordant to fibered knots, as
follows.  The Alexander polynomial of a fibered knot is monic.  Consider   a knot
$K$ with $\Delta_K(t) = 2t^2 - 3t +2$.  If $K$ were concordant to a fibered knot, then
$\Delta_K(t)g(t) = f(t)f(t^{-1})$ for some monic polynomial $g$ and integral $f$.  However,
since $\Delta_K(t)$ is irreducible and symmetric, it would have to be a factor of
$f(t)$ and of $f(t^{-1})$, giving it even exponent in
$\Delta_K(t)g(t)$, implying it is a factor of $g(t)$, contradicting monotonicity.

As mentioned above, Soma~\cite{so} proved that fibered knots are  concordant to hyperbolic
fibered knots.

The most significant result associating fibering and concordance is the theorem of Casson and
Gordon~\cite{cg3}.

\begin{theorem} If $K$ is a fibered ribbon knot, then the monodromy of the fibration extends
over some solid handlebody.  \end{theorem}

\subsection{Unknotting Number.} The unknotting number of a knot
$K$  is the least number of crossing changes that must be made in any diagram of $K$ to convert
it to an unknot.  This is closely related to the 4--ball genus  of a
knot (see the discussion above) and questions regarding the slicing of a knot in manifolds
bounded by $S^3$ other than
$B^4$, for instance, a once punctured connected sums of copies of $S^2 \times S^2$.  A related
invariant that is more closely tied to  concordance  was introduced by Askitas~\cite{as, o},
which we call the slicing number of a knot:
$u_s(K)$ is the minimum number of crossing changes required to convert a knot into a slice
knot.   It is relatively easy to see that the 4--ball genus of a knot provides a lower bound on
the slicing number; it was shown in~\cite{my1} and later in~\cite{l3} that these two need not be
equal.

\section{Problems}

 Past problem sets that include questions related to the knot concordance group
include~\cite{go2, k}. \vskip.2in

\begin{enumerate}

\item   {\bf Is every slice knot a ribbon knot?}

     A knot is ribbon if it bounds an embedded
disk in
$B^4$ having no local maxima (with respect to the radial function) in its interior.  In the
topological category this is not defined, so one asks the following instead: {\bf is every slice
knot    homotopically ribbon?} (That is, does
$K$ bound a disk $D$ in
$B^4$ such that
$\pi_1(S^3 - K) \to \pi_1(B^4 - D)$ is surjective?) In the smooth setting one then has the
additional question: {\bf is every homotopically ribbon knot a ribbon knot?} 

One has little basis to conjecture here.  Perhaps obstructions will arise (in either category)
but the lack of potential examples is discouraging.  On the other hand, topological surgery
might provide a proof in that category, but would give little indication concerning the smooth
setting.

\item  {\bf Describe all torsion in  $\calc$.} 

  Beginning with~\cite{fm} the question of whether there is any odd torsion has been open. More
generally, the only known torsion in
$\calc$ is two torsion that arises from knots that are concordant to negative amphicheiral knots,
and Conjecture~\ref{gconject} (first suggested in~\cite{go2}) states that  negative
amphicheirality   is the source of all (two) torsion in
$\calc$. 

As described in Section~\ref{3dprops} the Seifert form $$V_3 =
 \left(  \begin{tabular}{c   c} 
      1  &   1   \\ 
   0   &  $-3$  \\
\end{tabular}
\right) $$ represents 2--torsion in $\calg$ but cannot be represented by a negative amphicheiral
knot.   

The prospects for understanding 4--torsion look better.  A start has been made in~\cite{ ln2,
ln1} where it is shown, for instance, that no knot with Seifert form 
$$V_5 =
 \left(  \begin{tabular}{c   c} 
      1  &   1   \\ 
   0   &  $-5$  \\
\end{tabular}
\right) $$ can be of order 4 in $\calc$, although every such knot is of order 4 in
$\calg$.

Closely related to questions of torsion is the question: {\bf Does Levine's homomorphism split?}
That is, is there a homomorphism
$\psi \co \calg \to \calc$ such that $\phi \circ \psi   $ is the identity?  An affirmative answer
would yield elements of order 4 in
$\calc$ as well as elements of order 2 that do not arise from negative amphicheiral knots.

See~\cite{kw, mo} for computations of the algebraic orders of small crossing number knots.

\item {\bf If the knots $K$ and $K \# J$ are doubly slice, that is cross-sections of {\em
unknotted} 2--spheres in $ \rr^4$, is $J$ doubly slice?}  

The study of double knot concordance has a long history, with some of the initial work appearing
in~\cite{su}.  Other references include~\cite{ gl4, ke2, ke3, le4, le5, st1}.   The property of
double sliceness can be used to define a double concordance group  which maps onto
$\calc$ and there is a corresponding algebraic double concordance group formed using
quotienting by the set of hyperbolic Seifert forms rather than metabolic forms.  Algebraic
invariants  show that the kernel is infinitely generated, and Casson-Gordon invariants and
Cochran-Orr-Teichner methods apply in the case that algebraic invariants do 
not~\cite{gl4, kim2}.   Although a variety of questions regarding double null concordance can be
asked, this problem points to the underlying geometric difficulty of the topic.

\item {\bf Describe the structure of   the kernel of Levine's homomorphism,   $  \cala =
\mbox{ker}(\phi \co \calc \to
\calg)$.}

 It is known~\cite{ji, l5} that $\cala$ contains a subgroup isomorphic to
$\zz^\infty
\oplus
\zz_2^\infty$.  A reasonable conjecture is that $\cala \cong 
\zz^\infty \oplus
\zz_2^\infty$.  It has recently been shown by the author~\cite{l8}  that results of Ozsv\'ath and
Szab\'o~\cite{os} imply that $\cala$ has a summand isomorphic to $\zz$.  This implies that
$\cala$ contains elements that are not divisible and that $\cala$ is not a divisible group.
 There remains the unlikely possibility that $\cala$ does contain infinitely divisible elements,
perhaps including summands isomorphic to $\qq$ and $\qq / \zz$.

\item {\bf Describe the kernel of the map from $\calc$ to the topological concordance group,
$\calc_{top}$.}

  It is known that the kernel is nontrivial, containing for instance non-smoothly
slice Alexander polynomial one knots.  (See~\cite{cg, gom} for early references.) In    fact
it contains an infinitely generated such subgroup~\cite{e}. What more can be said about this
kernel?

\item{\bf Identify new relationships between the various unknotting numbers and genera of a
knot.} 

Here is a problem that seems to test the limits of presently known techniques.  If $K$
can be converted into a slice knot by making $m$ positive crossing changes and $n$ negative
crossing changes, then a geometric construction yields a surface bounded by $K$ in the 4--ball
of genus max$\{m,n\}$. Conversely: {\bf If the 4--ball genus of $K$ is $g_4$, can $K$ be
converted into a slice knot by making $g_4$ positive and $g_4$ negative crossing changes?} 
A simpler question  ask the same thing except for  the 3--sphere genus $g_3$ instead of $g_4$. 
(It is interesting to note that at this time it seems unknown if  the classical unknotting
number satisfies
$u(K)
\le 2 g_3(K)$.)
\end{enumerate}


 \newcommand{\etalchar}[1]{$^{#1}$}

\end{document}